\documentclass{amsart}%
\usepackage{amssymb}
\usepackage{amsfonts}
\usepackage{geometry}%
\usepackage{amsmath}%
\usepackage{graphicx}
\providecommand{\U}[1]{\protect\rule{.1in}{.1in}}
\newtheorem{theorem}{Theorem}
\theoremstyle{plain}
\newtheorem{acknowledgement}[theorem]{Acknowledgement}

\newtheorem{lemma}[theorem]{Lemma}

\newtheorem{proposition}[theorem]{Proposition}
\newtheorem{remark}[theorem]{Remark}

\numberwithin{equation}{section}
\begin{document}
\title[Martingale averages and Fourier extension]{Martingale averages, Fourier extension and the Schr\"{o}dinger equation}
\author{Lyudmila Korobenko}
\address{Reed College\\
Portland, Oregon, USA}
\email{lkoroben@math.mcmaster.ca}
\author{Cristian Rios}
\address{University of Calgary\\
Calgary, Alberta, Canada}
\email{crios@ucalgary.ca}
\author[E. T. Sawyer]{Eric T. Sawyer$^{\dagger}$}
\address{Eric T. Sawyer, Department of Mathematics and Statistics\\
McMaster University\\
1280 Main Street West\\
Hamilton, Ontario L8S 4K1 Canada}
\email{sawyer@mcmaster.ca}
\thanks{$\dagger$ Research supported in part by a grant from the National Science and
Engineering Research Council of Canada.}
\date{\today }

\begin{abstract}
The martingale averaged Fourier extension theorem holds for the paraboloid in
$n$ dimensions for $\left(  p,q\right)  \in A_{\operatorname*{Knapp}}$, i.e.%
\[
\mathbb{E}_{\pm}\left\Vert \mathcal{E}_{S}f_{\pm}\right\Vert _{L^{q}\left(
\mathbb{R}^{n}\right)  }\leq C\left\Vert f\right\Vert _{L^{p}\left(  U\right)
}\ ,\ \ \ \ \ \left(  p,q\right)  \in A_{\operatorname*{Knapp}},
\]
where%
\[
A_{\operatorname*{Knapp}}\equiv\left\{  \left(  p,q\right)  :q=\frac{n+1}%
{n-1}p^{\prime},2<p<\frac{2n}{n-1}\right\}  .
\]
and where $\mathcal{E}_{S}f_{\alpha}\left(  \xi\right)  \equiv\int
_{U}e^{-i\Phi\left(  x\right)  \cdot\xi}f_{\alpha}\left(  x\right)  dx$,
$\Phi\left(  x\right)  =\left(  x,\left\vert x\right\vert ^{2}\right)  $,
$f_{\pm}$ is the random martingale transform of the Haar expansion of $f$.
Using embedding and interpolation, this gives the complete range of possible
exponents (the diagonal case $p=q$ was obtained in \cite{Saw7} using smooth
Alpert wavelets instead). OpenAI has recently announced proofs of the
deterministic analogue on the diagonal $p=q$ \cite{Ope} and of the local
smoothing conjecture \cite{OpeLS}, in all dimensions. Since a deterministic
estimate at a point of $A_{\operatorname*{Knapp}}$ implies the probabilistic
one there, the point of the present paper is a short, elementary and
self-contained proof of the probabilistic estimate, which uses no Kakeya-type input.

We give an application of the probabilistic Knapp theorem to the global
initial value problem for the Schr\"{o}dinger equation, in which we obtain
sharp averaged $L^{p}\rightarrow L^{q}$ regularity results averaged over
martingale transforms of the initial data $f$.

\end{abstract}
\maketitle
\tableofcontents

\section{Introduction}

The Fourier extension conjecture in $n\geq2$ dimensions arose from unpublished
work of E. Stein in 1967, see e.g. \cite[see the Notes at the end of Chapter
IX, p. 432, where Stein proved the restriction conjecture for $1\leq
p<\frac{4d}{3d+1}$]{Ste2} and \cite{Ste},%
\begin{equation}
\left\Vert \widehat{f\sigma_{n-1}}\right\Vert _{L^{q}\left(  \mathbb{R}%
^{n}\right)  }\leq C_{p}\left\Vert f\right\Vert _{L^{p}\left(  \sigma
_{n-1}\right)  },\ \text{for }f\in L^{p}\left(  \sigma_{n-1}\right)
,\label{extension}%
\end{equation}
for%
\begin{align}
& \left(  \frac{1}{p},\frac{1}{q}\right)  \in\left[  0,1\right]  ^{2}\text{
lying strictly below the line }\frac{1}{q}=\frac{n-1}{2n}%
,\label{Knapp conditions}\\
& \text{and on or below the Knapp line through }\left(  \frac{n-1}{2n}%
,\frac{n-1}{2n}\right)  \text{ and }\left(  1,0\right)  ,\nonumber
\end{align}
and where $\sigma_{n-1}$ is surface measure on $S$ and where $S$ is either the
sphere $\mathbb{S}^{n-1}$ or a smooth compactly supported measure on the
paraboloid $\mathbb{P}^{n-1}$. The case $n=2$ of the Fourier extension
conjecture for the circle was proved over half a century ago by C. Fefferman
\cite{Fef}, see also A. Zygmund \cite{Zyg} and L. Carleson and P. Sj\"{o}lin
\cite{CaSj}. Just recently, OpenAI has announced a proof of the diagonal
inequality $p=q$ in all dimensions with $L\exists\forall N$ certification
\cite{Ope}, as well as a proof of the local smoothing conjecture in all
dimensions, with a Lean formalization, which is stated there to imply the
Kakeya conjecture and Stein's restriction conjecture \cite{OpeLS}. A proof of
the full Fourier extension conjecture for the sphere in dimensions $n\geq3$,
including the Knapp line, has been claimed in \cite{Chu}. The present paper
does not rely on any of these results.

We will assume $S$ is a smooth compact piece of the paraboloid in the
arguments in this paper, but in the end, the size of $S$ plays no role, so our
theorem below is stated globally. Using the inequalities of Khintchine and
Kahane, the probabilistic inequality (\ref{prob inequ}) can be equivalently
reformulated as a Fourier square function inequality,%
\begin{equation}
\left\Vert \left(  \sum_{I\in\mathcal{G}}\left\vert T\triangle_{I}f\right\vert
^{2}\right)  ^{\frac{1}{2}}\right\Vert _{L^{q}\left(  \mathbb{R}^{n}\right)
}\leq C\left\Vert f\right\Vert _{L^{p}\left(  \mathbb{R}^{n-1}\right)
},\label{square function}%
\end{equation}
where $T$ is the extension operator (\ref{def T par}) and $\triangle_{I}$ are
the Haar projections, see Subsection \ref{setting}. It was pointed out in
\cite{Saw7} that using the square function inequality (\ref{square function}),
the probabilistic analogue (\ref{prob inequ}) in Theorem \ref{main} below
fails for the same pairs $\left(  p,q\right)  $ that (\ref{extension}) is
currently known to fail for - namely those not satisfying
(\ref{Knapp conditions}). Before the announcements mentioned above, the Knapp
arc
\[
A_{\operatorname*{Knapp}}\equiv\left\{  \left(  p,q\right)  :\frac{1}{p}%
+\frac{n+1}{n-1}\frac{1}{q}=1\text{ and }2<p<\frac{2n}{n-1}\right\}
\]
was the only set of pairs $\left(  p,q\right)  $ for which it was \emph{not}
known whether (\ref{prob inequ}) holds. Note also that the Knapp arc can be
described as the set of points $\left(  p,q\right)  $ such that $\left(
\frac{1}{p},\frac{1}{q}\right)  $ lies on the line segment
$S_{\operatorname*{Knapp}}$ that is the subset of the line $x+\frac{n+1}%
{n-1}y=1$ lying strictly between the points corresponding to $p=2$ and
$p=\frac{2n}{n-1}$. The martingale averaged Fourier extension theorem below
will be used to obtain probabilistic regularity results for the global
Schr\"{o}dinger equation.

It should be emphasized that randomizing over all dyadic scales removes the
difficult Kakeya structure in the deterministic conjecture, and hence that the
theorem below does not deal with any Kakeya issues.

In this paper we work with the classical Haar multipliers%
\[
\mathcal{H}_{\mathbf{a}}^{\mathcal{G}}f\equiv\sum_{I\in\mathcal{G}}%
a_{I}\triangle_{I}f,\ \ \ \ \ \triangle_{I}f=\sum_{e\in\Gamma}\left\langle
f,h_{I}^{e}\right\rangle h_{I}^{e},
\]
see Subsection \ref{setting}, instead of the smooth Alpert wavelets used in
\cite{Saw7} to prove the corresponding probabilistic diagonal case. The Haar
multipliers are the unsmoothed case $\eta=0$, $\kappa=1$ of the Alpert
multipliers in \cite{Saw7}. The price paid is that a single Haar wave packet
decays only to first order away from its tube, and the proof compensates for
this by using the product structure of both the Haar functions and the phase
of the paraboloid. The smooth Alpert multipliers of \cite{Saw7} are treated by
the same argument in Remark \ref{smooth Alpert}.

\begin{theorem}
[probabilistic Knapp estimate for Haar multipliers]\label{main}Let $n\geq2$,
let $\mathcal{G}\subset\mathcal{D}$, and let $\left(  p,q\right)  \in
A_{\operatorname*{Knapp}}$. Then for $f\in L^{p}\cap L^{2}\left(
\mathbb{R}^{n-1}\right)  $,%
\begin{equation}
\mathbb{E}\left\Vert T\mathcal{H}_{\mathbf{a}}^{\mathcal{G}}f\right\Vert
_{L^{q}\left(  \mathbb{R}^{n}\right)  }\leq C\left\Vert f\right\Vert
_{L^{p}\left(  \mathbb{R}^{n-1}\right)  },\label{prob inequ}%
\end{equation}
with $C$ depending only on $n$ and $p$; in particular $C$ does not depend on
$\mathcal{G}$. For infinite $\mathcal{G}$, $T\mathcal{H}_{\mathbf{a}%
}^{\mathcal{G}}f$ is the limit of finite sums described in Subsection
\ref{proof main}.
\end{theorem}

We make two remarks on the statement.

\begin{enumerate}
\item The endpoint $p=2$, $q=\frac{2\left(  n+1\right)  }{n-1}$ is included as
well: the proof below applies verbatim there, and in any case
(\ref{prob inequ}) then follows from the Strichartz estimate \cite{Str}, since
$\left\Vert \mathcal{H}_{\mathbf{a}}^{\mathcal{G}}f\right\Vert _{L^{2}}%
\leq\left\Vert f\right\Vert _{L^{2}}$.

\item Any deterministic estimate on the Knapp line implies the corresponding
probabilistic estimate. Indeed, if $\left\Vert Tg\right\Vert _{L^{q}\left(
\mathbb{R}^{n}\right)  }\leq C\left\Vert g\right\Vert _{L^{p}\left(
\mathbb{R}^{n-1}\right)  }$ for some $\left(  p,q\right)  \in
A_{\operatorname*{Knapp}}$, then by Burkholder's inequality (H2) in Subsection
\ref{setting},%
\[
\mathbb{E}\left\Vert T\mathcal{H}_{\mathbf{a}}^{\mathcal{G}}f\right\Vert
_{L^{q}\left(  \mathbb{R}^{n}\right)  }\leq\sup_{\mathbf{a}}\left\Vert
T\mathcal{H}_{\mathbf{a}}^{\mathcal{G}}f\right\Vert _{L^{q}\left(
\mathbb{R}^{n}\right)  }\leq C\sup_{\mathbf{a}}\left\Vert \mathcal{H}%
_{\mathbf{a}}^{\mathcal{G}}f\right\Vert _{L^{p}\left(  \mathbb{R}%
^{n-1}\right)  }\leq C\,C_{n,p}\left\Vert f\right\Vert _{L^{p}\left(
\mathbb{R}^{n-1}\right)  }.
\]
So Theorem \ref{main} would also follow from a deterministic Knapp line
estimate for the paraboloid. The point of the present paper is a short,
elementary and self-contained proof of the probabilistic estimate, which uses
no Kakeya-type input: after randomization, only a bush of tubes through the
origin has to be controlled, see Subsection \ref{kakeya}.
\end{enumerate}

\begin{acknowledgement}
The current paper benefited from conversations with the AI agent Claude,
especially in two points,
\end{acknowledgement}

\begin{enumerate}
\item using familiar Haar wavelets instead of the smooth Alpert wavelets as in
\cite{Saw7}, and

\item starting the proof for the Knapp segment with the square function
formulation instead of the averaged formulation, resulting in a greatly
simplified and shortened argument as compared to that in \cite{Saw7} for the
diagonal case using smooth Alpert wavelets. However, the probabilistic
analogue for the \emph{sphere} may require smooth wavelets, as the phase is no
longer a sum of one-dimensional phases - see Subsubsection \ref{sphere} below.
\end{enumerate}

\section{The global initial value problem for the Schr\"{o}dinger equation}

Denote by $\mathcal{S}\left(  \mathbb{R}^{n}\right)  $ the space of rapidly
decreasing functions on $\mathbb{R}^{n}$, and let $\mathcal{S}\left(
\mathbb{R}^{n}\right)  ^{\ast}$ be the space of tempered distributions.
Consider the global initial value problem for the Schr\"{o}dinger equation,%
\begin{equation}
\left\{
\begin{array}
[c]{cccc}%
i\frac{\partial u}{\partial t}\left(  t,x\right)  -\bigtriangleup_{x}u\left(
t,x\right)  & = & 0 & \left(  t,x\right)  \in\mathbb{R}\times\mathbb{R}%
^{n-1}\\
u\left(  0,x\right)  & = & f\left(  x\right)  & x\in\mathbb{R}^{n-1}%
\end{array}
\right.  ,\label{IVP Schrodinger}%
\end{equation}
where $f\in\mathcal{S}\left(  \mathbb{R}^{n-1}\right)  ^{\ast}$ is a tempered
distribution on $\mathbb{R}^{n-1}$. We will now show the probabilistic Fourier
extension theorem can be used to prove a probabilistic regularity theorem for
(\ref{IVP Schrodinger}).

But first we recall the following well known connection between the
Schr\"{o}dinger equation and Fourier extensions. Let $\Phi\left(  \xi\right)
\equiv\left(  2\pi\left\vert \xi\right\vert ^{2},\xi\right)  $, $\xi
\in\mathbb{R}^{n-1}$, be the standard parameterization of $\mathbb{P}^{n-1}$
(with the time coordinate first), and let $\mathcal{F}_{d}$ denote the Fourier
transform in Euclidean space $\mathbb{R}^{d}$, for $d=n-1$ and $n$.

\begin{theorem}
Let $f\in\mathcal{S}\left(  \mathbb{R}^{n-1}\right)  ^{\ast}$ be a tempered
distribution in $\mathbb{R}^{n-1}$. Then $u_{f}\equiv\mathcal{F}_{n}^{-1}%
\Phi_{\ast}\mathcal{F}_{n-1}f$ is a tempered distribution and is a solution to
the Schr\"{o}dinger initial value problem (\ref{IVP Schrodinger}).
\end{theorem}

In this section we write%
\[
\mathcal{E}g\left(  t,x\right)  \equiv\mathcal{F}_{n}^{-1}\Phi_{\ast}g\left(
t,x\right)  =\int_{\mathbb{R}^{n-1}}e^{2\pi i\left(  2\pi t\left\vert
\xi\right\vert ^{2}+x\cdot\xi\right)  }g\left(  \xi\right)  d\xi
,\ \ \ \ \ \left(  t,x\right)  \in\mathbb{R}\times\mathbb{R}^{n-1},
\]
so that $u_{f}=\mathcal{EF}_{n-1}f$. This parameterization differs from the
one used in Section \ref{short} only by a linear change of variables: if $T$
is the extension operator in (\ref{def T par}), then $\mathcal{E}g\left(
t,x\right)  =Tg\left(  -x,-2\pi t\right)  $, and hence%
\[
\left\Vert \mathcal{E}g\right\Vert _{L^{q}\left(  \mathbb{R}^{n}\right)
}=\left(  2\pi\right)  ^{-\frac{1}{q}}\left\Vert Tg\right\Vert _{L^{q}\left(
\mathbb{R}^{n}\right)  }.
\]
In particular, Theorem \ref{main} applies to $\mathcal{E}$ with the same exponents.

Both estimates below concern the whole paraboloid, and the dilations
$\xi\mapsto\lambda\xi$ of $\mathbb{R}^{n-1}$ force the exponents onto the
Knapp line $q=\frac{n+1}{n-1}p^{\prime}$; off this line no such global
estimate can hold. For $p=2$ the relevant estimate is that of Strichartz
\cite{Str}, the global form of the Stein-Tomas theorem \cite{Tom}, see e.g.
\cite{Tao}. Together with Plancherel's theorem it gives the following $L^{2}$
regularity result for the initial value problem (\ref{IVP Schrodinger}),%
\[
\left\Vert u_{f}\right\Vert _{L^{q}\left(  \mathbb{R}^{n}\right)  }=\left\Vert
\mathcal{EF}_{n-1}f\right\Vert _{L^{q}\left(  \mathbb{R}^{n}\right)  }\leq
C_{n}\left\Vert \mathcal{F}_{n-1}f\right\Vert _{L^{2}\left(  \mathbb{R}%
^{n-1}\right)  }=C_{n}\left\Vert f\right\Vert _{L^{2}\left(  \mathbb{R}%
^{n-1}\right)  },\ \ \ \ \ q=\frac{2\left(  n+1\right)  }{n-1},
\]
for all $f\in L^{2}\left(  \mathbb{R}^{n-1}\right)  $.

Using the Hausdorff-Young inequality $\left\Vert \mathcal{F}_{n-1}f\right\Vert
_{L^{p}}\leq\left\Vert f\right\Vert _{L^{p^{\prime}}}$, $p\geq2$, instead, we
obtain from the martingale averaged estimate in Theorem \ref{main} the
following sharp probabilistic $L^{p^{\prime}}\rightarrow L^{q}$ regularity
estimate. For $\left(  p,q\right)  \in A_{\operatorname*{Knapp}}$ and $f\in
L^{p^{\prime}}\cap L^{2}\left(  \mathbb{R}^{n-1}\right)  $, so that
$\mathcal{F}_{n-1}f\in L^{p}\cap L^{2}\left(  \mathbb{R}^{n-1}\right)  $,%
\begin{align*}
\mathbb{E}\left\Vert u_{\mathcal{F}_{n-1}^{-1}\mathcal{H}_{\mathbf{a}%
}^{\mathcal{G}}\mathcal{F}_{n-1}f}\right\Vert _{L^{q}\left(  \mathbb{R}%
^{n}\right)  }  & =\mathbb{E}\left\Vert \mathcal{F}_{n}^{-1}\Phi_{\ast
}\mathcal{H}_{\mathbf{a}}^{\mathcal{G}}\mathcal{F}_{n-1}f\right\Vert
_{L^{q}\left(  \mathbb{R}^{n}\right)  }=\mathbb{E}\left\Vert \mathcal{EH}%
_{\mathbf{a}}^{\mathcal{G}}\mathcal{F}_{n-1}f\right\Vert _{L^{q}\left(
\mathbb{R}^{n}\right)  }\\
& \leq C_{n,p}\left\Vert \mathcal{F}_{n-1}f\right\Vert _{L^{p}\left(
\mathbb{R}^{n-1}\right)  }\leq C_{n,p}\left\Vert f\right\Vert _{L^{p^{\prime}%
}\left(  \mathbb{R}^{n-1}\right)  }.
\end{align*}

\section{Proof of the probabilistic Knapp estimate\label{short}}

In this section we prove the probabilistic estimate on the open Knapp segment
Theorem \ref{main} in full detail. Besides two classical facts about the Haar
system (its tensor product structure, and Burkholder's martingale
inequalities), the proof uses only Khintchine's inequality, Young's
convolution inequality, the Fefferman--Stein inequality, and van der Corput's
lemma for one-dimensional oscillatory integrals (Lemma \ref{chirp}). We
include many routine details, and after each step we comment in lay language
on what it does and why it is needed.

Throughout, $A\lesssim B$ means $A\leq CB$ with a constant $C$ depending only
on $n$ and $p$. The constant never depends on $f$, on the coefficients, on the
collection of cubes, on the scale $s$, or on the point $\xi$. We write
$A\approx B$ when $A\lesssim B\lesssim A$.

\subsection{Setting, notation and preliminaries\label{setting}}

We introduce the extension operator, cubes, scales and exponents we consider,
as well as the\ Haar system and multipliers, along with a standard toolkit of
inequalities, including Khintchine, Young, Fefferman-Stein and van der Corput.

\textbf{The extension operator.} We use the paraboloid $\Phi\left(  x\right)
=\left(  x,\left\vert x\right\vert ^{2}\right)  $, $x\in\mathbb{R}^{n-1}$, and
write points of $\mathbb{R}^{n}$ as $\xi=\left(  \xi^{\prime},\xi_{n}\right)
$ with $\xi^{\prime}\in\mathbb{R}^{n-1}$ and $\xi_{n}\in\mathbb{R}$. The
extension operator is%
\begin{equation}
Tf\left(  \xi\right)  =\int_{\mathbb{R}^{n-1}}e^{-2\pi i\Phi\left(  x\right)
\cdot\xi}f\left(  x\right)  dx=\int_{\mathbb{R}^{n-1}}e^{-2\pi i\left(
x\cdot\xi^{\prime}+\left\vert x\right\vert ^{2}\xi_{n}\right)  }f\left(
x\right)  dx.\label{def T par}%
\end{equation}
Thus $\xi_{n}$ is the direction normal to the paraboloid at the origin. The
sphere is discussed briefly in Subsubsection \ref{sphere}.

\textbf{Cubes and scales.} $\mathcal{D}$ is the standard dyadic grid of
$\mathbb{R}^{n-1}$. For $I\in\mathcal{D}$ we write $\ell\left(  I\right)  $
for its side length, $x_{I}$ for its center, and $\left\vert I\right\vert
=\ell\left(  I\right)  ^{n-1}$ for its volume. Each $I\in\mathcal{D}$ is a
product $I=I_{1}\times\cdots\times I_{n-1}$ of dyadic intervals of length
$\ell\left(  I\right)  $. We fix an arbitrary subcollection $\mathcal{G}%
\subset\mathcal{D}$. One may take $\mathcal{G}=\mathcal{G}\left[  S\right]  $,
the dyadic subcubes of a fixed cube $S$, but nothing below depends on this
choice, and $\mathcal{G}=\mathcal{D}$ is allowed. For $s\in\mathbb{Z}$ we
write $\mathcal{D}_{s}=\left\{  I\in\mathcal{D}:\ell\left(  I\right)
=2^{-s}\right\}  $ for the cubes of scale $s$. Larger $s$ means smaller cubes.

\textbf{Exponents.} Fix $\left(  p,q\right)  \in A_{\operatorname*{Knapp}}$,
that is, $2<p<n^{\ast}=\frac{2n}{n-1}$ and $q=\frac{n+1}{n-1}p^{\prime}$. We
use three elementary consequences.

\begin{enumerate}
\item[(E1)] $2\leq p<q$. Indeed, $q>p$ if and only if $\frac{n+1}{n-1}%
>\frac{p}{p^{\prime}}=p-1$, that is, $p<n^{\ast}$. In particular $\frac{p}%
{2}\geq1$ and $\frac{q}{2}\geq1$.

\item[(E2)] $0<\frac{2}{p}-\frac{2}{q}<1$. Positivity is (E1). Since $\frac
{2}{q}=\frac{2\left(  n-1\right)  }{n+1}\left(  1-\frac{1}{p}\right)  $, the
difference is an increasing function of $\frac{1}{p}$, so it is largest at
$p=2$, where it equals $1-\frac{n-1}{n+1}=\frac{2}{n+1}<1$. We put%
\[
\gamma\equiv\left(  n-1\right)  \left(  \frac{2}{p}-\frac{2}{q}\right)  >0.
\]

\item[(E3)] The Knapp identity: $\frac{n-1}{2}\left(  1-\frac{q}{p^{\prime}%
}\right)  =-1$, because $\frac{q}{p^{\prime}}=\frac{n+1}{n-1}$.
\end{enumerate}

\textbf{The Haar system and Haar multipliers.} On $\mathbb{R}$, for a dyadic
interval $J$ with dyadic halves $J_{-}$ and $J_{+}$, put%
\[
h_{J}^{0}\equiv\left\vert J\right\vert ^{-\frac{1}{2}}\mathbf{1}%
_{J}\ \ \ \ \ \text{and}\ \ \ \ \ h_{J}^{1}\equiv\left\vert J\right\vert
^{-\frac{1}{2}}\left(  \mathbf{1}_{J_{-}}-\mathbf{1}_{J_{+}}\right)  .
\]
To keep letters distinct we write $e\in\Gamma\equiv\left\{  0,1\right\}
^{n-1}\setminus\left\{  0\right\}  $ for the index set of Haar functions, and
we reserve $\mathbf{a}=\left\{  a_{I}\right\}  $ for the random signs. For
$I=I_{1}\times\cdots\times I_{n-1}\in\mathcal{D}$ and $e\in\Gamma$ put%
\[
h_{I}^{e}\left(  x\right)  \equiv\prod_{i=1}^{n-1}h_{I_{i}}^{e_{i}}\left(
x_{i}\right)  ,\ \ \ \ \ \triangle_{I}^{e}f\equiv\left\langle f,h_{I}%
^{e}\right\rangle h_{I}^{e},\ \ \ \ \ \triangle_{I}f\equiv\sum_{e\in\Gamma
}\triangle_{I}^{e}f,
\]
and for $\mathbf{a}=\left\{  a_{I}\right\}  \in\left\{  -1,1\right\}
^{\mathcal{G}}$ the Haar multiplier is $\mathcal{H}_{\mathbf{a}}^{\mathcal{G}%
}f=\sum_{I\in\mathcal{G}}a_{I}\triangle_{I}f$. The signs $a_{I}$ are
independent and equal to $\pm1$ with probability $\frac{1}{2}$ each (the
product measure $\mu$), and $\mathbb{E}$ denotes expectation. We use two
classical facts.

\begin{enumerate}
\item[(H1)] (Tensor structure.) $\left\{  h_{I}^{e}\right\}  _{I\in
\mathcal{D},e\in\Gamma}$ is an orthonormal basis of $L^{2}\left(
\mathbb{R}^{n-1}\right)  $. Each $h_{I}^{e}$ is a product of one-dimensional
functions $h_{I_{i}}^{e_{i}}$, each of which is $\ell\left(  I\right)
^{-\frac{1}{2}}$ times a signed sum of the indicators of at most two intervals
contained in $I_{i}$. Moreover $\left\vert h_{I}^{e}\right\vert =\left\vert
I\right\vert ^{-\frac{1}{2}}\mathbf{1}_{I}$.

\item[(H2)] (Burkholder.) For $s\in\mathbb{Z}$, $\sum_{I\in\mathcal{D}_{s}%
}\triangle_{I}f=\mathbb{E}_{s+1}f-\mathbb{E}_{s}f$ is the martingale
difference of the dyadic filtration. Hence for $1<p<\infty$ and $f\in
L^{p}\left(  \mathbb{R}^{n-1}\right)  $,%
\[
\sup_{\varepsilon\in\left\{  -1,1\right\}  ^{\mathcal{D}}}\left\Vert
\sum_{I\in\mathcal{D}}\varepsilon_{I}\triangle_{I}f\right\Vert _{L^{p}%
}+\left\Vert \left(  \sum_{I\in\mathcal{D}}\left\vert \triangle_{I}%
f\right\vert ^{2}\right)  ^{\frac{1}{2}}\right\Vert _{L^{p}}\leq
C_{n,p}\left\Vert f\right\Vert _{L^{p}},
\]
and the Haar system is an unconditional basis of $L^{p}\left(  \mathbb{R}%
^{n-1}\right)  $. In particular $\mathcal{H}_{\mathbf{a}}^{\mathcal{G}}$ is
bounded on $L^{p}$ uniformly in $\mathbf{a}$ and $\mathcal{G}$.
\end{enumerate}

\textbf{Standard tools.}

\begin{enumerate}
\item[(T1)] (Khintchine's inequality.) For $0<q<\infty$ there is $B_{q}$ such
that for all complex numbers $w_{1},\dots,w_{M}$ and independent random signs
$a_{1},\dots,a_{M}$,%
\[
\left(  \mathbb{E}\left\vert \sum_{i=1}^{M}a_{i}w_{i}\right\vert ^{q}\right)
^{\frac{1}{q}}\leq B_{q}\left(  \sum_{i=1}^{M}\left\vert w_{i}\right\vert
^{2}\right)  ^{\frac{1}{2}}.
\]
(See e.g.\ Grafakos, \emph{Classical Fourier Analysis}, Appendix C. The
complex case follows from the real one by treating real and imaginary parts separately.)

\item[(T2)] (Young's inequality.) If $1\leq t,u,v\leq\infty$ and $1+\frac
{1}{v}=\frac{1}{t}+\frac{1}{u}$, then $\left\Vert K\ast m\right\Vert _{L^{v}%
}\leq\left\Vert K\right\Vert _{L^{t}}\left\Vert m\right\Vert _{L^{u}} $, on
$\mathbb{R}^{n-1}$ and equally for sequences on $\mathbb{Z}$.

\item[(T3)] (Fefferman--Stein.) If $M$ is the Hardy--Littlewood maximal
operator and $1<p<\infty$, then%
\[
\left\Vert \left(  \sum_{k}\left(  Mg_{k}\right)  ^{2}\right)  ^{\frac{1}{2}%
}\right\Vert _{L^{p}}\leq C\left\Vert \left(  \sum_{k}\left\vert
g_{k}\right\vert ^{2}\right)  ^{\frac{1}{2}}\right\Vert _{L^{p}}.
\]

\item[(T4)] (Two elementary inequalities.) (i) If $x_{k}\geq0$ and $r\geq1$,
then $\sum_{k}x_{k}^{r}\leq\left(  \sum_{k}x_{k}\right)  ^{r}$. To see this,
put $X=\sum_{k}x_{k}$; then $\sum_{k}\left(  x_{k}/X\right)  ^{r}\leq\sum
_{k}x_{k}/X=1$. (ii) If $\beta_{k}\geq0$ and $p\leq q$, then $\left(  \sum
_{k}\beta_{k}^{q}\right)  ^{\frac{1}{q}}\leq\left(  \sum_{k}\beta_{k}%
^{p}\right)  ^{\frac{1}{p}}$. This is (i) with $x_{k}=\beta_{k}^{p}$ and
$r=\frac{q}{p}$.

\item[(T5)] (van der Corput's lemma, \cite[Chapter VIII, Proposition 2]%
{Ste2}.) If $\phi$ is real on $\left[  \alpha,\beta\right]  $ and $\left\vert
\phi^{\prime\prime}\right\vert \geq\lambda>0$ there, then $\left\vert
\int_{\alpha}^{\beta}e^{i\phi}\right\vert \leq8\lambda^{-\frac{1}{2}}$.
\end{enumerate}

\subsection{The picture behind the proof\label{picture}}

Fix a cube $I$ of side $\ell$ and center $x_{I}$. Lemma \ref{packet} below
says that, up to tails of first order in each coordinate, $\left\vert
Th_{I}^{e}\right\vert $ is a \emph{wave packet}:

\begin{itemize}
\item For $\left\vert \xi_{n}\right\vert \lesssim\ell^{-2}$ it has size at
most $\approx\left\vert I\right\vert ^{\frac{1}{2}}$, and comparable to
$\left\vert I\right\vert ^{\frac{1}{2}}$ on a fixed fraction of the tube%
\[
\left\{  \xi:\left\vert \xi_{n}\right\vert \lesssim\ell^{-2},\ \left\vert
\xi^{\prime}+2\xi_{n}x_{I}\right\vert \lesssim\ell^{-1}\right\}  .
\]
This is a cylinder of radius $\ell^{-1}$ and length $\ell^{-2}$ around the
line through the origin in direction $\left(  -2x_{I},1\right)  $. That
direction is normal to the paraboloid at $\Phi\left(  x_{I}\right)  $, so this
is the familiar Knapp tube. (The lower bound is not proved here; only the
upper bounds of Lemma \ref{packet} are used below.)

\item For $\left\vert \xi_{n}\right\vert \gtrsim\ell^{-2}$, that is, beyond
the \emph{coherence length} $\ell^{-2}$, the packet disperses. It fills the
cone $\left\{  \left\vert \xi^{\prime}+2\xi_{n}x_{I}\right\vert \lesssim
\left\vert \xi_{n}\right\vert \ell\right\}  $, of aperture $\approx\ell$, with
amplitude at most $\approx\left\vert I\right\vert ^{\frac{1}{2}}\left(
\ell^{2}\left\vert \xi_{n}\right\vert \right)  ^{-\frac{n-1}{2}}$.
\end{itemize}

As a sanity check, on each horizontal slice $\left\{  \xi_{n}=\mathrm{const}%
\right\}  $ we have (size)$^{2}\times$(area)$\ \approx\left\vert I\right\vert
\rho^{-\left(  n-1\right)  }\left(  \rho/\ell\right)  ^{n-1}=1$, where
$\rho=1+\ell^{2}\left\vert \xi_{n}\right\vert $. This agrees with Plancherel's
theorem in $\xi^{\prime}$: the function $\xi^{\prime}\mapsto Th\left(
\xi^{\prime},\xi_{n}\right)  $ is the Fourier transform of $x\mapsto e^{-2\pi
i\left\vert x\right\vert ^{2}\xi_{n}}h\left(  x\right)  $, whose $L^{2}$ norm
is $\left\Vert h\right\Vert _{L^{2}}=1$.

Two features of this picture drive the proof.

\begin{enumerate}
\item \emph{Every packet passes through the origin.} This is because
$h_{I}^{e}$ is not modulated. A modulated function $e^{2\pi ix\cdot w}h\left(
x\right)  $ would produce a packet translated by $w$. In a Haar expansion such
modulations are not carried by single Haar functions: they are built from the
coefficients of many cubes, and the random signs scramble them.

\item \emph{The random signs remove interference.} The packets $G_{I}$ coming
from different cubes interfere in $Tf=\sum_{I}G_{I}$, so $Tf$ can be much
larger or much smaller than the individual packets. In the deterministic
conjecture, controlling this interference is the whole difficulty; it is where
Kakeya-type problems enter. By Khintchine's inequality, the average over signs
of $\left\Vert \sum_{I}a_{I}G_{I}\right\Vert _{L^{q}}$ is controlled by
$\left\Vert \left(  \sum_{I}\left\vert G_{I}\right\vert ^{2}\right)
^{\frac{1}{2}}\right\Vert _{L^{q}}$, which only sees the sizes $\left\vert
G_{I}\right\vert $. After this step the problem is \emph{positive}: we must
estimate the $L^{q/2}$ norm of the nonnegative function $F=\sum_{I}\left\vert
G_{I}\right\vert ^{2}$, a weighted sum of bumps adapted to a bush of tubes
through the origin.
\end{enumerate}

\textbf{Natural coordinates.} For $\xi_{n}\neq0$ set%
\[
z=-\frac{\xi^{\prime}}{2\xi_{n}}\in\mathbb{R}^{n-1}.
\]
The phase $x\cdot\xi^{\prime}+\left\vert x\right\vert ^{2}\xi_{n}$ is
stationary in $x$ exactly at $x=z$. Equivalently, $z$ is the point of the
parameter space whose normal line to the paraboloid passes through $\xi$. In
these coordinates the packet of a cube $I$ of scale $s$ is visible at $\xi$
only if $\mathrm{dist}\left(  z,I\right)  \lesssim r_{s}\left(  \xi
_{n}\right)  $, where%
\[
r_{s}\left(  \xi_{n}\right)  =\frac{1}{2}\left(  \frac{2^{s}}{\left\vert
\xi_{n}\right\vert }+2^{-s}\right)  .
\]
The first term is the tube radius $\ell^{-1}=2^{s}$ divided by the height
$\left\vert \xi_{n}\right\vert $, i.e.\ the angular size of the tube as seen
from height $\xi_{n}$. The second term is the aperture of the dispersive cone.
The minimum over $s$ is about $\left\vert \xi_{n}\right\vert ^{-\frac{1}{2}}$.
It is attained at the \emph{resonant scale} $2^{-s}\approx\left\vert \xi
_{n}\right\vert ^{-\frac{1}{2}}$, for which the height $\left\vert \xi
_{n}\right\vert $ equals the coherence length $\ell^{-2}$. Writing
$\sigma\left(  \xi_{n}\right)  =\frac{1}{2}\log_{2}\left\vert \xi
_{n}\right\vert $ for the resonant scale, we have $r_{s}\left(  \xi
_{n}\right)  \approx\left\vert \xi_{n}\right\vert ^{-\frac{1}{2}}2^{\left\vert
s-\sigma\left(  \xi_{n}\right)  \right\vert }$.

\subsection{The six step proof}

\begin{enumerate}
\item Randomization and Khintchine's inequality reduce the proof to an $L^{q}
$ norm estimate of a square function.

\item Fix a height $\xi_{n}$. Apart from the Jacobian factor $\left(
2\left\vert \xi_{n}\right\vert \right)  ^{-\left(  n-1\right)  }$, the
function $F$ restricted to that height is a sum over scales $s$ of the
coefficient densities $m_{s}$, each blurred at radius $r_{s}\left(  \xi
_{n}\right)  $.

\item Blurring at radius $r$ improves integrability from $L^{p/2}$ to
$L^{q/2}$, with a gain $r^{-\gamma}$, by Young's inequality. The farther $s$
is from the resonant scale, the larger the blur and therefore the gain, which
is $2^{-\gamma\left\vert s-\sigma\left(  \xi_{n}\right)  \right\vert } $.

\item The Knapp relation makes the remaining power of $\left\vert \xi
_{n}\right\vert $ exactly $-1$, so what is left is the dilation-invariant
measure $d\xi_{n}/\left\vert \xi_{n}\right\vert $.

\item Summing over dyadic heights then becomes a discrete convolution in the
scale parameter, which Young's inequality on $\mathbb{Z}$ controls.

\item Finally, since $q\geq p\geq2$, the $\ell^{q}$ sum over scales is bounded
by the $L^{p}$ norm of the square function.
\end{enumerate}

\subsubsection{Wave packets through the origin}

We first estimate the one-dimensional integrals into which $Th_{I}^{e}$ factors.

\begin{lemma}
[one-dimensional chirp integrals]\label{chirp}For real $v,\mu$ and $\left[
\alpha,\beta\right]  \subset\left[  -\frac{1}{2},\frac{1}{2}\right]  $,%
\begin{equation}
j\left(  v,\mu\right)  \equiv\left\vert \int_{\alpha}^{\beta}e^{-i\left(
vu+\mu u^{2}\right)  }du\right\vert \leq27\left(  1+\left\vert \mu\right\vert
\right)  ^{-\frac{1}{2}}\left(  1+\frac{\left\vert v\right\vert }{1+\left\vert
\mu\right\vert }\right)  ^{-1}.\label{chirp bound}%
\end{equation}

\end{lemma}

\begin{proof}
Let $\phi\left(  u\right)  =vu+\mu u^{2}$, so that $\phi^{\prime}\left(
u\right)  =v+2\mu u$ and $\phi^{\prime\prime}=2\mu$.

\emph{(a) The dispersive bound.} Trivially $j\leq1$. If $\mu\neq0$, van der
Corput's lemma (T5) with $\left\vert \phi^{\prime\prime}\right\vert
=2\left\vert \mu\right\vert $ gives $j\leq8\left(  2\left\vert \mu\right\vert
\right)  ^{-\frac{1}{2}}<6\left\vert \mu\right\vert ^{-\frac{1}{2}}$. Hence%
\begin{equation}
j\left(  v,\mu\right)  \leq\min\left\{  1,6\left\vert \mu\right\vert
^{-\frac{1}{2}}\right\}  \leq9\left(  1+\left\vert \mu\right\vert \right)
^{-\frac{1}{2}}\ \ \ \ \ \text{for all }v\in\mathbb{R}.\label{dispersive 1d}%
\end{equation}

\emph{(b) Decay away from the stationary point.} Suppose $\left\vert
v\right\vert \geq2\left\vert \mu\right\vert $ and $v\neq0$. For $\left\vert
u\right\vert \leq\frac{1}{2}$ we have $\left\vert 2\mu u\right\vert
\leq\left\vert \mu\right\vert \leq\frac{1}{2}\left\vert v\right\vert $, so
$\left\vert \phi^{\prime}\right\vert \geq\frac{1}{2}\left\vert v\right\vert $
on $\left[  \alpha,\beta\right]  $, $\phi^{\prime}$ has constant sign there,
and $1/\phi^{\prime}$ is monotone. Integrating by parts,%
\[
\int_{\alpha}^{\beta}e^{-i\phi}du=\left[  \frac{e^{-i\phi}}{-i\phi^{\prime}%
}\right]  _{\alpha}^{\beta}-\int_{\alpha}^{\beta}e^{-i\phi}\left(  \frac
{1}{-i\phi^{\prime}}\right)  ^{\prime}du,
\]
so that%
\begin{equation}
j\left(  v,\mu\right)  \leq\frac{4}{\left\vert v\right\vert }+\left\vert
\frac{1}{\phi^{\prime}\left(  \beta\right)  }-\frac{1}{\phi^{\prime}\left(
\alpha\right)  }\right\vert \leq\frac{6}{\left\vert v\right\vert
}\ \ \ \ \ \text{when }\left\vert v\right\vert \geq2\left\vert \mu\right\vert
.\label{nonstat 1d}%
\end{equation}
This first order decay is sharp: it comes from the endpoints $\alpha,\beta$,
i.e.\ from the jumps of the Haar functions.

\emph{(c) Conclusion.} Put $\Lambda=1+\left\vert \mu\right\vert $. If
$\left\vert v\right\vert \leq2\Lambda$, then $1+\left\vert v\right\vert
/\Lambda\leq3$, and (\ref{dispersive 1d}) gives $j\leq27\Lambda^{-\frac{1}{2}%
}\left(  1+\left\vert v\right\vert /\Lambda\right)  ^{-1}$. If $\left\vert
v\right\vert >2\Lambda$, then $\left\vert v\right\vert \geq2\left\vert
\mu\right\vert $, so (\ref{nonstat 1d}) gives $j\leq6/\left\vert v\right\vert
$, while $1+\left\vert v\right\vert /\Lambda\leq2\left\vert v\right\vert
/\Lambda$ gives%
\[
\Lambda^{-\frac{1}{2}}\left(  1+\frac{\left\vert v\right\vert }{\Lambda
}\right)  ^{-1}\geq\frac{\Lambda^{\frac{1}{2}}}{2\left\vert v\right\vert }%
\geq\frac{1}{2\left\vert v\right\vert }.
\]
Hence $j\leq12\Lambda^{-\frac{1}{2}}\left(  1+\left\vert v\right\vert
/\Lambda\right)  ^{-1}$ in this case too.
\end{proof}

\begin{lemma}
[Haar wave packets through the origin]\label{packet}Let $I\in\mathcal{D}$ have
side $\ell$ and center $x_{I}$, let $e\in\Gamma$, and put $\rho\equiv
1+\ell^{2}\left\vert \xi_{n}\right\vert $. Then for all $\xi=\left(
\xi^{\prime},\xi_{n}\right)  \in\mathbb{R}^{n}$,%
\begin{equation}
\left\vert Th_{I}^{e}\left(  \xi\right)  \right\vert \leq54^{n-1}\left\vert
I\right\vert ^{\frac{1}{2}}\rho^{-\frac{n-1}{2}}\prod_{i=1}^{n-1}\left(
1+\frac{\ell\left\vert \xi_{i}^{\prime}+2\xi_{n}\left(  x_{I}\right)
_{i}\right\vert }{\rho}\right)  ^{-1}.\label{packet bound}%
\end{equation}

\end{lemma}

Each factor $\left(  1+\cdot\right)  ^{-1}$ is $\approx1$ exactly when
$\left\vert \xi_{i}^{\prime}+2\xi_{n}\left(  x_{I}\right)  _{i}\right\vert
\lesssim\rho/\ell=\ell^{-1}+\ell\left\vert \xi_{n}\right\vert $. This is the
tube for $\left\vert \xi_{n}\right\vert \leq\ell^{-2}$ and the cone for
$\left\vert \xi_{n}\right\vert \geq\ell^{-2}$, as in Subsection \ref{picture}.
The factor $\rho^{-\frac{n-1}{2}}$ is the dispersive decay. Away from the tube
or cone, the packet decays only to first order in each coordinate. This cannot
be improved, because $h_{I}^{e}$ has jumps, and it is the reason why the
kernel in Step 2 below is a product kernel.

\begin{proof}
Since $x\cdot\xi^{\prime}+\left\vert x\right\vert ^{2}\xi_{n}=\sum_{i=1}%
^{n-1}\left(  x_{i}\xi_{i}^{\prime}+x_{i}^{2}\xi_{n}\right)  $, the tensor
structure (H1) gives%
\[
Th_{I}^{e}\left(  \xi\right)  =\prod_{i=1}^{n-1}\int_{\mathbb{R}}e^{-2\pi
i\left(  t\xi_{i}^{\prime}+t^{2}\xi_{n}\right)  }h_{I_{i}}^{e_{i}}\left(
t\right)  dt.
\]
Fix $i$ and write $c=\left(  x_{I}\right)  _{i}$. By (H1) the $i$-th factor is
$\ell^{-\frac{1}{2}}$ times a signed sum of at most two integrals of $e^{-2\pi
i\left(  t\xi_{i}^{\prime}+t^{2}\xi_{n}\right)  }$ over intervals contained in
$I_{i}$. Substitute $t=c+\ell u$, with $u$ in an interval $\left[
\alpha,\beta\right]  \subset\left[  -\frac{1}{2},\frac{1}{2}\right]  $, and
expand%
\[
t\xi_{i}^{\prime}+t^{2}\xi_{n}=\left(  c\xi_{i}^{\prime}+c^{2}\xi_{n}\right)
+\ell u\left(  \xi_{i}^{\prime}+2\xi_{n}c\right)  +\ell^{2}\xi_{n}u^{2}.
\]
Each such integral therefore has modulus $\ell\,j\left(  v,\mu\right)  $,
where%
\begin{equation}
v=2\pi\ell\left(  \xi_{i}^{\prime}+2\xi_{n}c\right)  ,\ \ \ \ \ \mu=2\pi
\ell^{2}\xi_{n}.\label{def v mu}%
\end{equation}
In words, $j$ is the integral of an indicator against a linear frequency $v$
and a quadratic \textquotedblleft chirp\textquotedblright\ $\mu$. The
frequency $v$ measures the distance from $\xi$ to the axis of the tube in
units of $\ell^{-1}$; the chirp $\mu$ measures the height in units of the
coherence length $\ell^{-2}$. Note that%
\begin{equation}
\rho\leq1+\left\vert \mu\right\vert \leq2\pi\rho,\ \ \ \ \ \text{so}%
\ \ \ \ \ \left(  1+\left\vert \mu\right\vert \right)  ^{-\frac{1}{2}}\leq
\rho^{-\frac{1}{2}}\ \ \ \ \ \text{and}\ \ \ \ \ \frac{\left\vert v\right\vert
}{1+\left\vert \mu\right\vert }\geq\frac{\ell\left\vert \xi_{i}^{\prime}%
+2\xi_{n}c\right\vert }{\rho}.\label{mu rho}%
\end{equation}
By Lemma \ref{chirp} and (\ref{mu rho}), the $i$-th factor is at most%
\[
\ell^{-\frac{1}{2}}\cdot2\cdot\ell\cdot27\rho^{-\frac{1}{2}}\left(
1+\frac{\ell\left\vert \xi_{i}^{\prime}+2\xi_{n}c\right\vert }{\rho}\right)
^{-1}=54\,\ell^{\frac{1}{2}}\rho^{-\frac{1}{2}}\left(  1+\frac{\ell\left\vert
\xi_{i}^{\prime}+2\xi_{n}c\right\vert }{\rho}\right)  ^{-1}.
\]
Taking the product over $i$ and using $\ell^{\frac{n-1}{2}}=\left\vert
I\right\vert ^{\frac{1}{2}}$ gives (\ref{packet bound}).
\end{proof}

\subsubsection{The key estimate}

It is convenient to separate the Haar machinery from the harmonic analysis.
For a finite collection $\mathcal{G}_{0}\subset\mathcal{D}$ and complex
coefficients $c=\left\{  c_{I,e}\right\}  _{I\in\mathcal{G}_{0},e\in\Gamma}$,
define%
\begin{align}
G_{I}  & \equiv\sum_{e\in\Gamma}c_{I,e}Th_{I}^{e}=T\left(  \sum_{e\in\Gamma
}c_{I,e}h_{I}^{e}\right)  ,\label{def m_s}\\
m_{s}\left(  x\right)   & \equiv\sum_{I\in\mathcal{G}_{0}\cap\mathcal{D}_{s}%
}\sum_{e\in\Gamma}\frac{\left\vert c_{I,e}\right\vert ^{2}}{\left\vert
I\right\vert }\mathbf{1}_{I}\left(  x\right)  ,\ \ \ \ \ \mathcal{S}%
c\equiv\left(  \sum_{s\in\mathbb{Z}}m_{s}\right)  ^{\frac{1}{2}}.\nonumber
\end{align}
We call $m_{s}$ the \emph{coefficient density at scale }$s$, and
$\mathcal{S}c$ the \emph{coefficient square function}.

\begin{proposition}
[key estimate]\label{key}Let $n\geq2$ and $\left(  p,q\right)  \in
A_{\operatorname*{Knapp}}$. There is $C$, depending only on $n$ and $p$, such
that for every finite $\mathcal{G}_{0}\subset\mathcal{D}$ and all coefficients
$c$,%
\[
\mathbb{E}\left\Vert \sum_{I\in\mathcal{G}_{0}}a_{I}G_{I}\right\Vert
_{L^{q}\left(  \mathbb{R}^{n}\right)  }\leq C\left\Vert \mathcal{S}%
c\right\Vert _{L^{p}\left(  \mathbb{R}^{n-1}\right)  }.
\]

\end{proposition}

\begin{lemma}
[coefficient square function]\label{coef sq}For $1<p<\infty$ and $g\in
L^{p}\cap L^{2}\left(  \mathbb{R}^{n-1}\right)  $,%
\[
\left\Vert \left(  \sum_{I\in\mathcal{D}}\sum_{e\in\Gamma}\frac{\left\vert
\left\langle g,h_{I}^{e}\right\rangle \right\vert ^{2}}{\left\vert
I\right\vert }\mathbf{1}_{I}\right)  ^{\frac{1}{2}}\right\Vert _{L^{p}}\leq
C\left\Vert g\right\Vert _{L^{p}}.
\]

\end{lemma}

Theorem \ref{main} follows by applying Proposition \ref{key} to the
coefficients $c_{I,e}=\left\langle f,h_{I}^{e}\right\rangle $, for which
$\sum_{I}a_{I}G_{I}=T\mathcal{H}_{\mathbf{a}}^{\mathcal{G}}f$, and then
bounding $\left\Vert \mathcal{S}c\right\Vert _{L^{p}}$ by Lemma \ref{coef sq}.
The details, including the passage to infinite $\mathcal{G}$, are in
Subsection \ref{proof main}.

Proposition \ref{key} is a statement about arbitrary coefficient sequences.
Haar functions enter it only through the size and decay of single packets
(Lemma \ref{packet}). No cancellation between different cubes is used, and
neither is the vanishing mean of $h_{I}^{e}$. The heart of the matter is
Proposition \ref{key}.

\subsubsection{Proof of Proposition \ref{key}}

Fix a finite $\mathcal{G}_{0}\subset\mathcal{D}$ and coefficients $c$; all
sums below are finite.

\medskip

\textbf{Step 1 (randomization: from signs to a positive quantity).} Let%
\[
F\left(  \xi\right)  \equiv\sum_{I\in\mathcal{G}_{0}}\sum_{e\in\Gamma
}\left\vert c_{I,e}\right\vert ^{2}\left\vert Th_{I}^{e}\left(  \xi\right)
\right\vert ^{2}.
\]
We claim that%
\begin{equation}
\mathbb{E}\left\Vert \sum_{I\in\mathcal{G}_{0}}a_{I}G_{I}\right\Vert _{L^{q}%
}\leq B_{q}\left\vert \Gamma\right\vert ^{\frac{1}{2}}\left\Vert F^{\frac
{1}{2}}\right\Vert _{L^{q}},\ \ \ \ \ \left\vert \Gamma\right\vert
=2^{n-1}-1.\label{step1}%
\end{equation}
By H\"{o}lder's inequality in the probability space, $\mathbb{E}\left\Vert
X\right\Vert _{L^{q}}\leq\left(  \mathbb{E}\left\Vert X\right\Vert _{L^{q}%
}^{q}\right)  ^{\frac{1}{q}}$. By Tonelli's theorem,%
\[
\mathbb{E}\left\Vert \sum_{I}a_{I}G_{I}\right\Vert _{L^{q}}^{q}=\int
_{\mathbb{R}^{n}}\mathbb{E}\left\vert \sum_{I}a_{I}G_{I}\left(  \xi\right)
\right\vert ^{q}d\xi\leq B_{q}^{q}\int_{\mathbb{R}^{n}}\left(  \sum
_{I}\left\vert G_{I}\left(  \xi\right)  \right\vert ^{2}\right)  ^{\frac{q}%
{2}}d\xi,
\]
where we applied Khintchine's inequality (T1) at each fixed $\xi$ with
$w_{I}=G_{I}\left(  \xi\right)  $. Finally, by Cauchy--Schwarz over the
$\left\vert \Gamma\right\vert $ terms in $G_{I}$,%
\[
\left\vert G_{I}\left(  \xi\right)  \right\vert ^{2}=\left\vert \sum
_{e\in\Gamma}c_{I,e}Th_{I}^{e}\left(  \xi\right)  \right\vert ^{2}%
\leq\left\vert \Gamma\right\vert \sum_{e\in\Gamma}\left\vert c_{I,e}%
\right\vert ^{2}\left\vert Th_{I}^{e}\left(  \xi\right)  \right\vert ^{2}.
\]
Combining these three facts gives (\ref{step1}).

\emph{Commentary.} This is the only place where the random signs are used.
Nothing essential is lost: Kahane's inequality and the lower Khintchine
inequality give $\mathbb{E}\left\Vert \sum_{I}a_{I}G_{I}\right\Vert _{L^{q}%
}\approx\left\Vert \left(  \sum_{I}\left\vert G_{I}\right\vert ^{2}\right)
^{\frac{1}{2}}\right\Vert _{L^{q}}$, so the probabilistic estimate is
\emph{equivalent} to a square function estimate, and $F$ is (up to the
Cauchy--Schwarz in $e$) the square of that square function. From now on
everything is positive, and only upper bounds for $\left\vert Th_{I}%
^{e}\right\vert $ are used.

\medskip

\textbf{Step 2 (pointwise bound: a blurred picture of the coefficients).} For
$s\in\mathbb{Z}$, $\xi_{n}\neq0$, $r>0$ and $z\in\mathbb{R}^{n-1}$, define%
\[
\rho_{s}\left(  \xi_{n}\right)  \equiv1+4^{-s}\left\vert \xi_{n}\right\vert
,\ \ \ \ \ r_{s}\left(  \xi_{n}\right)  \equiv\frac{\rho_{s}\left(  \xi
_{n}\right)  }{2^{1-s}\left\vert \xi_{n}\right\vert }=\frac{1}{2}\left(
\frac{2^{s}}{\left\vert \xi_{n}\right\vert }+2^{-s}\right)  ,
\]
and the \emph{product kernel}%
\begin{equation}
K_{r}\left(  z\right)  \equiv r^{-\left(  n-1\right)  }\prod_{i=1}%
^{n-1}\left(  1+\frac{\left\vert z_{i}\right\vert }{r}\right)  ^{-2}%
.\label{def K}%
\end{equation}
Thus $\rho_{s}$ is the $\rho$ of Lemma \ref{packet} for cubes of scale $s$,
and $K_{r}\left(  z\right)  =r^{-\left(  n-1\right)  }K_{1}\left(  z/r\right)
$ is a bump of width $r$. Since each factor $\left(  1+\left\vert
z_{i}\right\vert \right)  ^{-2}$ is integrable on $\mathbb{R}$, we have for
every $t\geq1$,%
\begin{equation}
\left\Vert K_{1}\right\Vert _{L^{t}}=\left(  \int_{\mathbb{R}}\left(
1+\left\vert y\right\vert \right)  ^{-2t}dy\right)  ^{\frac{n-1}{t}}=\left(
\frac{2}{2t-1}\right)  ^{\frac{n-1}{t}}<\infty,\label{K norm}%
\end{equation}
and in particular $\left\Vert K_{r}\right\Vert _{L^{1}}=\left\Vert
K_{1}\right\Vert _{L^{1}}=2^{n-1}$. We claim that for $\xi_{n}\neq0$ and
$z\in\mathbb{R}^{n-1}$,%
\begin{equation}
F\left(  -2\xi_{n}z,\xi_{n}\right)  \lesssim\left(  2\left\vert \xi
_{n}\right\vert \right)  ^{-\left(  n-1\right)  }\sum_{s\in\mathbb{Z}}\left(
K_{r_{s}\left(  \xi_{n}\right)  }\ast m_{s}\right)  \left(  z\right)
.\label{step2}%
\end{equation}

Fix $\xi_{n}\neq0$ and $z$, and put $\xi=\left(  -2\xi_{n}z,\xi_{n}\right)  $.
Fix also $s$ and a cube $I\in\mathcal{G}_{0}\cap\mathcal{D}_{s}$, so that
$\ell=\ell\left(  I\right)  =2^{-s}$, and abbreviate $\rho=\rho_{s}\left(
\xi_{n}\right)  $ and $r=r_{s}\left(  \xi_{n}\right)  $. Since $\xi
_{i}^{\prime}+2\xi_{n}\left(  x_{I}\right)  _{i}=-2\xi_{n}\left(
z_{i}-\left(  x_{I}\right)  _{i}\right)  $, we have for each $i$%
\[
\frac{\ell\left\vert \xi_{i}^{\prime}+2\xi_{n}\left(  x_{I}\right)
_{i}\right\vert }{\rho}=\frac{2\left\vert \xi_{n}\right\vert \ell}{\rho
}\left\vert z_{i}-\left(  x_{I}\right)  _{i}\right\vert =\frac{\left\vert
z_{i}-\left(  x_{I}\right)  _{i}\right\vert }{r},
\]
and Lemma \ref{packet} (squared) gives%
\begin{equation}
\left\vert Th_{I}^{e}\left(  \xi\right)  \right\vert ^{2}\lesssim\left\vert
I\right\vert \rho^{-\left(  n-1\right)  }\prod_{i=1}^{n-1}\left(
1+\frac{\left\vert z_{i}-\left(  x_{I}\right)  _{i}\right\vert }{r}\right)
^{-2}.\label{step2 i}%
\end{equation}
Next we replace the center $x_{I}$ by an arbitrary point $x\in I$, one
coordinate at a time. Since $r=\frac{1}{2\left\vert \xi_{n}\right\vert \ell
}+\frac{\ell}{2}\geq\frac{\ell}{2}\geq\left\vert x_{i}-\left(  x_{I}\right)
_{i}\right\vert $,%
\[
1+\frac{\left\vert z_{i}-x_{i}\right\vert }{r}\leq2+\frac{\left\vert
z_{i}-\left(  x_{I}\right)  _{i}\right\vert }{r}\leq2\left(  1+\frac
{\left\vert z_{i}-\left(  x_{I}\right)  _{i}\right\vert }{r}\right)  .
\]
Hence the product in (\ref{step2 i}) is at most $4^{n-1}\prod_{i}\left(
1+\left\vert z_{i}-x_{i}\right\vert /r\right)  ^{-2}$ for every $x\in I$, and
averaging over $x\in I$,%
\begin{equation}
\prod_{i=1}^{n-1}\left(  1+\frac{\left\vert z_{i}-\left(  x_{I}\right)
_{i}\right\vert }{r}\right)  ^{-2}\leq\frac{4^{n-1}}{\left\vert I\right\vert
}\int_{I}\prod_{i=1}^{n-1}\left(  1+\frac{\left\vert z_{i}-x_{i}\right\vert
}{r}\right)  ^{-2}dx.\label{step2 ii}%
\end{equation}
Multiply (\ref{step2 i}) by $\left\vert c_{I,e}\right\vert ^{2}$, use
(\ref{step2 ii}), and sum over $e\in\Gamma$ and $I\in\mathcal{G}_{0}%
\cap\mathcal{D}_{s}$. All these cubes have the same volume $\left\vert
I\right\vert =2^{-s\left(  n-1\right)  }$, so by the definition of $m_{s}$,%
\begin{align*}
\sum_{I\in\mathcal{G}_{0}\cap\mathcal{D}_{s}}\sum_{e\in\Gamma}\left\vert
c_{I,e}\right\vert ^{2}\left\vert Th_{I}^{e}\left(  \xi\right)  \right\vert
^{2}  & \lesssim2^{-s\left(  n-1\right)  }\rho^{-\left(  n-1\right)  }%
\int_{\mathbb{R}^{n-1}}m_{s}\left(  x\right)  \prod_{i=1}^{n-1}\left(
1+\frac{\left\vert z_{i}-x_{i}\right\vert }{r}\right)  ^{-2}dx\\
& =\left(  \frac{2^{-s}r}{\rho}\right)  ^{n-1}\left(  K_{r}\ast m_{s}\right)
\left(  z\right)  =\left(  2\left\vert \xi_{n}\right\vert \right)  ^{-\left(
n-1\right)  }\left(  K_{r}\ast m_{s}\right)  \left(  z\right)  ,
\end{align*}
where in the last step we used $2^{-s}r/\rho=1/\left(  2\left\vert \xi
_{n}\right\vert \right)  $, which is the definition of $r$. Summing over $s$
gives (\ref{step2}).

\emph{Commentary.} (a) The meaning of (\ref{step2}): at the point $\xi$, whose
stationary point is $z$, the scale-$s$ contribution is the coefficient density
$m_{s}$ averaged over a cube of side $\approx r_{s}\left(  \xi_{n}\right)  $
around $z$, with first order tails in each coordinate, times $\left(
2\left\vert \xi_{n}\right\vert \right)  ^{-\left(  n-1\right)  }$. That factor
is the Jacobian of $z\mapsto\xi^{\prime}=-2\xi_{n}z$; it spreads the
coefficient mass over the slice at height $\xi_{n}$.

(b) Step 2 loses nothing on average. By Plancherel in $\xi^{\prime}$ (see
Subsection \ref{picture}), for every $\xi_{n}$,%
\[
\int_{\mathbb{R}^{n-1}}F\left(  \xi^{\prime},\xi_{n}\right)  d\xi^{\prime
}=\sum_{I,e}\left\vert c_{I,e}\right\vert ^{2}\left\Vert h_{I}^{e}\right\Vert
_{L^{2}}^{2}=\sum_{I,e}\left\vert c_{I,e}\right\vert ^{2},
\]
while integrating the right side of (\ref{step2}) in $\xi^{\prime}$ gives%
\[
\left\Vert K_{1}\right\Vert _{L^{1}}\sum_{s}\int m_{s}=2^{n-1}\sum
_{I,e}\left\vert c_{I,e}\right\vert ^{2}.
\]
So (\ref{step2}) is sharp in the $L^{1}$ sense, and all the work is in
$L^{q/2}$.

(c) The product structure is essential. The decay in Lemma \ref{packet} is
only of first order in each coordinate, so the best \emph{radial} kernel
available would be $r^{-\left(  n-1\right)  }\left(  1+\left\vert z\right\vert
/r\right)  ^{-2}$, which is not integrable on $\mathbb{R}^{n-1}$ when $n\geq3
$. The product kernel (\ref{def K}) is in every $L^{t}$, $t\geq1$, because
each one-variable factor is integrable on $\mathbb{R}$.

(d) The two regimes of $r_{s}$: $r_{s}\approx2^{s}/\left\vert \xi
_{n}\right\vert $ when $4^{s}\geq\left\vert \xi_{n}\right\vert $ (we are
inside the tubes of scale $s$), and $r_{s}\approx2^{-s}$ when $4^{s}%
\leq\left\vert \xi_{n}\right\vert $ (we are in the dispersive cones). Writing
$\left\vert \xi_{n}\right\vert =4^{\sigma}$, i.e.\ $\sigma=\sigma\left(
\xi_{n}\right)  =\frac{1}{2}\log_{2}\left\vert \xi_{n}\right\vert $, we get
$2^{s}/\left\vert \xi_{n}\right\vert =2^{-\sigma}2^{s-\sigma}$ and
$2^{-s}=2^{-\sigma}2^{\sigma-s}$, hence%
\begin{equation}
r_{s}\left(  \xi_{n}\right)  =\frac{1}{2}\left\vert \xi_{n}\right\vert
^{-\frac{1}{2}}\left(  2^{s-\sigma}+2^{\sigma-s}\right)  \geq\frac{1}%
{2}\left\vert \xi_{n}\right\vert ^{-\frac{1}{2}}2^{\left\vert s-\sigma\left(
\xi_{n}\right)  \right\vert }.\label{r lower}%
\end{equation}

\medskip

\textbf{Step 3 (integration over horizontal slices and Young's inequality).}
By Fubini's theorem and the change of variables $\xi^{\prime}=-2\xi_{n}z$ at
fixed $\xi_{n}\neq0$, whose Jacobian is $\left(  2\left\vert \xi
_{n}\right\vert \right)  ^{n-1}$,%
\[
\left\Vert F^{\frac{1}{2}}\right\Vert _{L^{q}}^{q}=\int_{\mathbb{R}}%
\int_{\mathbb{R}^{n-1}}F\left(  \xi^{\prime},\xi_{n}\right)  ^{\frac{q}{2}%
}d\xi^{\prime}d\xi_{n}=\int_{\mathbb{R}}\left(  2\left\vert \xi_{n}\right\vert
\right)  ^{n-1}\int_{\mathbb{R}^{n-1}}F\left(  -2\xi_{n}z,\xi_{n}\right)
^{\frac{q}{2}}dz\,d\xi_{n}.
\]
Inserting (\ref{step2}),%
\begin{equation}
\left\Vert F^{\frac{1}{2}}\right\Vert _{L^{q}}^{q}\lesssim\int_{\mathbb{R}%
}\left(  2\left\vert \xi_{n}\right\vert \right)  ^{\left(  n-1\right)  \left(
1-\frac{q}{2}\right)  }\left\Vert \sum_{s\in\mathbb{Z}}K_{r_{s}\left(  \xi
_{n}\right)  }\ast m_{s}\right\Vert _{L^{q/2}}^{\frac{q}{2}}d\xi
_{n}.\label{step3a}%
\end{equation}
Now define $t$ by $1+\frac{2}{q}=\frac{1}{t}+\frac{2}{p}$, that is, $\frac
{1}{t}=1-\left(  \frac{2}{p}-\frac{2}{q}\right)  $. By (E2), $\frac{1}{t}%
\in\left(  0,1\right)  $, so $1<t<\infty$, and Young's inequality (T2) with
the exponents $\left(  t,\frac{p}{2},\frac{q}{2}\right)  $ (all $\geq1$ by
(E1)) gives%
\[
\left\Vert K_{r}\ast m_{s}\right\Vert _{L^{q/2}}\leq\left\Vert K_{r}%
\right\Vert _{L^{t}}\left\Vert m_{s}\right\Vert _{L^{p/2}}.
\]
Since $K_{r}\left(  z\right)  =r^{-\left(  n-1\right)  }K_{1}\left(
z/r\right)  $, a change of variables gives%
\[
\left\Vert K_{r}\right\Vert _{L^{t}}=r^{-\left(  n-1\right)  \left(
1-\frac{1}{t}\right)  }\left\Vert K_{1}\right\Vert _{L^{t}}=r^{-\gamma
}\left\Vert K_{1}\right\Vert _{L^{t}},
\]
and $\left\Vert K_{1}\right\Vert _{L^{t}}<\infty$ by (\ref{K norm}). Put
$b_{s}\equiv\left\Vert m_{s}\right\Vert _{L^{p/2}}$. By Minkowski's inequality
in $L^{q/2}$ (recall $\frac{q}{2}\geq1$) and (\ref{r lower}),%
\begin{equation}
\left\Vert \sum_{s\in\mathbb{Z}}K_{r_{s}\left(  \xi_{n}\right)  }\ast
m_{s}\right\Vert _{L^{q/2}}\leq\left\Vert K_{1}\right\Vert _{L^{t}}\sum
_{s\in\mathbb{Z}}r_{s}\left(  \xi_{n}\right)  ^{-\gamma}b_{s}\lesssim
\left\vert \xi_{n}\right\vert ^{\frac{\gamma}{2}}\sum_{s\in\mathbb{Z}%
}2^{-\gamma\left\vert s-\sigma\left(  \xi_{n}\right)  \right\vert }%
b_{s}.\label{step3b}%
\end{equation}

\emph{Commentary.} This is where the gain appears, and it is the only place
where $p<q$ is used in an essential way. An average over a cube of side $r$
maps $L^{p/2}$ to $L^{q/2}$ with norm $\approx r^{-\gamma}$: a function
concentrated on a set much smaller than the cube gets diluted. Scales far from
resonance are blurred over larger cubes (by (\ref{r lower})), so they gain
geometrically. At the endpoint $p=q=n^{\ast}$ we would have $\gamma=0$ and no
gain at all.

\medskip

\textbf{Step 4 (the Knapp relation: a scale-invariant measure).} Insert
(\ref{step3b}) into (\ref{step3a}). The power of $\left\vert \xi
_{n}\right\vert $ is%
\[
\left(  n-1\right)  \left(  1-\frac{q}{2}\right)  +\frac{\gamma q}{4}=\left(
n-1\right)  \left(  1-\frac{q}{2}+\frac{q}{2p}-\frac{1}{2}\right)  =\frac
{n-1}{2}\left(  1-q\left(  1-\frac{1}{p}\right)  \right)  =\frac{n-1}%
{2}\left(  1-\frac{q}{p^{\prime}}\right)  =-1,
\]
where we used $\frac{\gamma q}{4}=\left(  n-1\right)  \left(  \frac{q}%
{2p}-\frac{1}{2}\right)  $ and, in the last step, the Knapp identity (E3).
Hence%
\begin{equation}
\left\Vert F^{\frac{1}{2}}\right\Vert _{L^{q}}^{q}\lesssim\int_{\mathbb{R}%
}\left(  \sum_{s\in\mathbb{Z}}2^{-\gamma\left\vert s-\sigma\left(  \xi
_{n}\right)  \right\vert }b_{s}\right)  ^{\frac{q}{2}}\frac{d\xi_{n}%
}{\left\vert \xi_{n}\right\vert }.\label{step4}%
\end{equation}

\emph{Commentary.} The Knapp relation is used here and nowhere else. The
measure $d\xi_{n}/\left\vert \xi_{n}\right\vert $ is invariant under dilations
$\xi_{n}\mapsto\lambda\xi_{n}$. Rescaling the parameter space by $x\mapsto2x$
turns cubes of scale $s$ into cubes of scale $s-1$, divides the relevant
heights by $4$ (since $T\left[  f\left(  \cdot/2\right)  \right]  \left(
\xi^{\prime},\xi_{n}\right)  =2^{n-1}Tf\left(  2\xi^{\prime},4\xi_{n}\right)
$), and so lowers $\sigma$ by $1$. The weights $2^{-\gamma\left\vert
s-\sigma\right\vert }$ and the measure $d\xi_{n}/\left\vert \xi_{n}\right\vert
$ are unchanged, so (\ref{step4}) treats all scales alike. This is the scale
invariance of the Knapp line, and it is why no per-scale decay is available:
scales can only be summed through the off-diagonal decay $2^{-\gamma\left\vert
s-\sigma\right\vert }$.

\medskip

\textbf{Step 5 (summing over dyadic heights: a discrete convolution).} Split
$\mathbb{R}\setminus\left\{  0\right\}  =\bigcup_{j\in\mathbb{Z}}E_{j}$ with
$E_{j}=\left\{  4^{j}\leq\left\vert \xi_{n}\right\vert <4^{j+1}\right\}  $. On
$E_{j}$ we have $\sigma\left(  \xi_{n}\right)  \in\left[  j,j+1\right)  $, so
$\left\vert s-\sigma\right\vert \geq\left\vert s-j\right\vert -1$ and
$2^{-\gamma\left\vert s-\sigma\right\vert }\leq2^{\gamma}2^{-\gamma\left\vert
s-j\right\vert }$. Also $\int_{E_{j}}d\xi_{n}/\left\vert \xi_{n}\right\vert
=2\ln4$. Therefore%
\[
\left\Vert F^{\frac{1}{2}}\right\Vert _{L^{q}}^{q}\lesssim\sum_{j\in
\mathbb{Z}}\left(  \sum_{s\in\mathbb{Z}}k\left(  j-s\right)  b_{s}\right)
^{\frac{q}{2}},\ \ \ \ \ k\left(  i\right)  \equiv2^{-\gamma\left\vert
i\right\vert }.
\]
The inner sum is the discrete convolution $\left(  k\ast b\right)  \left(
j\right)  $. By Young's inequality on $\mathbb{Z}$ (recall $\frac{q}{2}\geq
1$), $\left\Vert k\ast b\right\Vert _{\ell^{q/2}}\leq\left\Vert k\right\Vert
_{\ell^{1}}\left\Vert b\right\Vert _{\ell^{q/2}}$, where $\left\Vert
k\right\Vert _{\ell^{1}}=\frac{1+2^{-\gamma}}{1-2^{-\gamma}}<\infty$ because
$\gamma>0$. Hence%
\begin{equation}
\left\Vert F^{\frac{1}{2}}\right\Vert _{L^{q}}^{q}\lesssim\sum_{s\in
\mathbb{Z}}b_{s}^{\frac{q}{2}}.\label{step5}%
\end{equation}

\emph{Commentary.} Each dyadic height $E_{j}$ is matched with its resonant
scale $s=j$, and scale $s$ contributes to height $j$ with weight
$2^{-\gamma\left\vert s-j\right\vert }$. The constant $\left\Vert k\right\Vert
_{\ell^{1}}\approx\frac{2}{\gamma\ln2}$ blows up as $\gamma\rightarrow0$,
i.e.\ as $p\rightarrow n^{\ast}$. This is the only source of the blow-up of
the constant in Theorem \ref{main}.

\medskip

\textbf{Step 6 (back to a single square function).} Write $\beta_{s}%
\equiv\left\Vert m_{s}^{1/2}\right\Vert _{L^{p}}$, so that $b_{s}^{\frac{q}%
{2}}=\left(  \int m_{s}^{\frac{p}{2}}\right)  ^{\frac{q}{p}}=\beta_{s}^{q}$.
Since $q\geq p$, (T4)(ii) gives $\sum_{s}\beta_{s}^{q}\leq\left(  \sum
_{s}\beta_{s}^{p}\right)  ^{\frac{q}{p}}$. Since $\frac{p}{2}\geq1$, (T4)(i)
gives the pointwise inequality $\sum_{s}m_{s}^{\frac{p}{2}}\leq\left(
\sum_{s}m_{s}\right)  ^{\frac{p}{2}}=\left(  \mathcal{S}c\right)  ^{p}$, and
integrating, $\sum_{s}\beta_{s}^{p}\leq\left\Vert \mathcal{S}c\right\Vert
_{L^{p}}^{p}$. Together with (\ref{step5}),%
\[
\left\Vert F^{\frac{1}{2}}\right\Vert _{L^{q}}\lesssim\left\Vert
\mathcal{S}c\right\Vert _{L^{p}},
\]
and Proposition \ref{key} follows from (\ref{step1}).

\emph{Commentary.} Both inequalities in Step 6 use $q\geq p\geq2$. They say
that, for $p\geq2$, an $\ell^{q}$ sum over scales of $L^{p}$ norms is
controlled by the $L^{p}$ norm of the $\ell^{2}$ sum over scales. For $p<2$
the second inequality fails.

\subsection{Proof of Lemma \ref{coef sq} and of Theorem \ref{main}%
\label{proof main}}

\begin{proof}
[Proof of Lemma \ref{coef sq}]Let $M$ be the Hardy--Littlewood maximal
operator over balls. Fix $I\in\mathcal{D}$ and $x\in I$. By orthonormality
(H1),%
\[
\sum_{e\in\Gamma}\frac{\left\vert \left\langle g,h_{I}^{e}\right\rangle
\right\vert ^{2}}{\left\vert I\right\vert }=\frac{1}{\left\vert I\right\vert
}\left\Vert \triangle_{I}g\right\Vert _{L^{2}}^{2}.
\]
The function $\triangle_{I}g$ is constant on each of the $2^{n-1}$ children
$J$ of $I$, with value $c_{J}$ there, say. Hence%
\[
\frac{1}{\left\vert I\right\vert }\left\Vert \triangle_{I}g\right\Vert
_{L^{2}}^{2}=2^{-\left(  n-1\right)  }\sum_{J}\left\vert c_{J}\right\vert
^{2}\leq2^{-\left(  n-1\right)  }\left(  \sum_{J}\left\vert c_{J}\right\vert
\right)  ^{2}=2^{n-1}\left(  \frac{1}{\left\vert I\right\vert }\int
_{I}\left\vert \triangle_{I}g\right\vert \right)  ^{2}.
\]
The ball $B=B\left(  x,\sqrt{n-1}\ell\left(  I\right)  \right)  $ contains $I$
and has volume $C_{n}\left\vert I\right\vert $, so $\frac{1}{\left\vert
I\right\vert }\int_{I}\left\vert \triangle_{I}g\right\vert \leq C_{n}M\left(
\triangle_{I}g\right)  \left(  x\right)  $. Hence, pointwise,%
\[
\sum_{I\in\mathcal{D}}\sum_{e\in\Gamma}\frac{\left\vert \left\langle
g,h_{I}^{e}\right\rangle \right\vert ^{2}}{\left\vert I\right\vert }%
\mathbf{1}_{I}\leq2^{n-1}C_{n}^{2}\sum_{I\in\mathcal{D}}\left(  M\triangle
_{I}g\right)  ^{2},
\]
and the lemma follows from the Fefferman--Stein inequality (T3) and the square
function bound in (H2).
\end{proof}

\emph{Commentary.} Since $\left\vert h_{I}^{e}\right\vert ^{2}=\left\vert
I\right\vert ^{-1}\mathbf{1}_{I}$, the left side of Lemma \ref{coef sq} is
exactly $\left\Vert \left(  \sum_{I,e}\left\vert \triangle_{I}^{e}g\right\vert
^{2}\right)  ^{\frac{1}{2}}\right\Vert _{L^{p}}$. The maximal function is used
only to pass from the individual pieces $\triangle_{I}^{e}g$ to the full
martingale differences $\triangle_{I}g$, to which Burkholder's inequality
applies: $\left\vert \triangle_{I}g\right\vert $ may be small on some children
of $I$, and the maximal function spreads its mass over all of $I$.

\begin{proof}
[Proof of Theorem \ref{main}]Put $c_{I,e}=\left\langle f,h_{I}^{e}%
\right\rangle $. For any collection $\mathcal{G}^{\prime}\subset\mathcal{D}$
write $\mathcal{S}_{\mathcal{G}^{\prime}}c=\left(  \sum_{I\in\mathcal{G}%
^{\prime}}\sum_{e\in\Gamma}\left\vert c_{I,e}\right\vert ^{2}\left\vert
I\right\vert ^{-1}\mathbf{1}_{I}\right)  ^{\frac{1}{2}}$.

\emph{Finite }$\mathcal{G}$\emph{.} For finite $\mathcal{G}_{0}\subset
\mathcal{G}$, each $h_{I}^{e}$ lies in $L^{1}$, so $\sum_{I\in\mathcal{G}_{0}%
}a_{I}G_{I}=\sum_{I\in\mathcal{G}_{0}}a_{I}T\triangle_{I}f=T\mathcal{H}%
_{\mathbf{a}}^{\mathcal{G}_{0}}f$. By Proposition \ref{key} and Lemma
\ref{coef sq},%
\begin{equation}
\mathbb{E}\left\Vert T\mathcal{H}_{\mathbf{a}}^{\mathcal{G}_{0}}f\right\Vert
_{L^{q}}\lesssim\left\Vert \mathcal{S}_{\mathcal{G}_{0}}c\right\Vert _{L^{p}%
}\leq\left\Vert \mathcal{S}_{\mathcal{D}}c\right\Vert _{L^{p}}\lesssim
\left\Vert f\right\Vert _{L^{p}}.\label{finite G}%
\end{equation}

\emph{Infinite }$\mathcal{G}$\emph{.} Choose finite collections $\mathcal{G}%
_{1}\subset\mathcal{G}_{2}\subset\cdots$ with union $\mathcal{G}$, and put
$X_{N}\left(  \mathbf{a}\right)  =T\mathcal{H}_{\mathbf{a}}^{\mathcal{G}_{N}%
}f$. For $N<N^{\prime}$, the difference $X_{N^{\prime}}-X_{N}$ is the same
expression for the finite collection $\mathcal{G}_{N^{\prime}}\setminus
\mathcal{G}_{N}$, so Proposition \ref{key} gives%
\[
\mathbb{E}\left\Vert X_{N^{\prime}}-X_{N}\right\Vert _{L^{q}}\lesssim
\left\Vert \mathcal{S}_{\mathcal{G}\setminus\mathcal{G}_{N}}c\right\Vert
_{L^{p}}\longrightarrow0\ \ \ \ \ \text{as }N\rightarrow\infty,
\]
by dominated convergence: $\mathcal{S}_{\mathcal{G}\setminus\mathcal{G}_{N}%
}c\leq\mathcal{S}_{\mathcal{D}}c\in L^{p}$, and $\mathcal{S}_{\mathcal{G}%
\setminus\mathcal{G}_{N}}c\rightarrow0$ pointwise as the tail of a convergent
series. So $\left(  X_{N}\right)  $ is a Cauchy sequence in the Bochner space
$L^{1}\left(  d\mu\left(  \mathbf{a}\right)  ;L^{q}\left(  \mathbb{R}%
^{n}\right)  \right)  $. Its limit $X$ satisfies $\mathbb{E}\left\Vert
X\right\Vert _{L^{q}}=\lim_{N}\mathbb{E}\left\Vert X_{N}\right\Vert _{L^{q}%
}\lesssim\left\Vert f\right\Vert _{L^{p}}$ by (\ref{finite G}), and it does
not depend on the exhaustion (interleave two exhaustions). We define
$T\mathcal{H}_{\mathbf{a}}^{\mathcal{G}}f=X\left(  \mathbf{a}\right)  $. When
$\mathcal{G}\subset\mathcal{G}\left[  S\right]  $ with $S$ bounded, this
agrees with the direct definition for $\mu$-almost every $\mathbf{a}$. Indeed,
the partial sums $\mathcal{H}_{\mathbf{a}}^{\mathcal{G}_{N}}f$ are supported
in $\overline{S}$ and, since the Haar system is an unconditional basis of
$L^{p}$ by (H2), they converge to $\mathcal{H}_{\mathbf{a}}^{\mathcal{G}}f$ in
$L^{p}$, hence in $L^{1}$, for every $\mathbf{a}$. Since $T$ maps $L^{1}$ to
$L^{\infty}$, $X_{N}\left(  \mathbf{a}\right)  \rightarrow T\mathcal{H}%
_{\mathbf{a}}^{\mathcal{G}}f$ uniformly, while a subsequence of $X_{N}\left(
\mathbf{a}\right)  $ converges to $X\left(  \mathbf{a}\right)  $ in $L^{q}$
for almost every $\mathbf{a}$.

Nothing in Lemma \ref{packet} or Proposition \ref{key} above uses $\ell\left(
I\right)  \leq1$, or that the cubes lie in a bounded set: $s$ ranges over all
of $\mathbb{Z}$. Hence Theorem \ref{main} holds with $\mathcal{G}=\mathcal{D}$
for the extension operator of the entire paraboloid. This is of course the
setting of the Schr\"{o}dinger application, which applies the theorem to
$\widehat{f}$ and uses Hausdorff-Young.
\end{proof}

\section{Concluding comments}

\subsection{Where the hypotheses are used\label{hyp}\ }

\begin{remark}
\begin{itemize}
\item The random signs are used only in Step 1.

\item $p\geq2$ is used in Young's inequality in Step 3 (it needs $\frac{p}%
{2}\geq1$) and in the superadditivity in Step 6.

\item $q>p$, i.e.\ $p<n^{\ast}$ is used through $\gamma>0$ in Steps 3 and 5.

\item $q\geq2$ is used in Minkowski's inequality in $L^{q/2}$ (Step 3) and in
Young's inequality on $\mathbb{Z}$ (Step 5).

\item The Knapp relation $q=\frac{n+1}{n-1}p^{\prime}$ is used only in Step 4,
through the exponent $-1$.

\item The tensor structure of the Haar functions, together with the splitting
of the phase of the paraboloid into one-variable phases is used only in Lemma
\ref{packet}, and through it in the product kernel (\ref{def K}) of Step 2.

\item The Haar structure in Theorem \ref{main} is used only through (H1) and
(H2), and this is the only place where the vanishing means of the $h_{I}^{e}$
are used.
\end{itemize}
\end{remark}

\subsection{Smooth Alpert multipliers\label{smooth Alpert}}

The same proof gives Theorem \ref{main} for the smooth Alpert multipliers
$\mathcal{A}_{\mathbf{a}}^{S_{\kappa,\eta}}f=\sum_{I\in\mathcal{G}}%
a_{I}\triangle_{I}^{\eta}f$ of \cite{Saw7}, where $\triangle_{I}^{\eta}%
f=\sum_{e}\left\langle S_{\kappa,\eta}^{-1}f,h_{I}^{e}\right\rangle
h_{I}^{e,\eta}$, $h_{I}^{e,\eta}=h_{I}^{e}\ast\phi_{\eta\ell\left(  I\right)
}$ are the smooth Alpert wavelets, and $S_{\kappa,\eta}$ is bounded and
invertible on $L^{p}$ for $\eta$ small (Theorem 6 of \cite{Saw7}), with
constants now depending also on $\kappa$ and $\eta$. Only Lemma \ref{packet}
and Lemma \ref{coef sq} change. Since $h_{I}^{e,\eta}$ is supported in
$\left(  1+\eta\right)  I$ and $\left\Vert \partial^{\alpha}h_{I}^{e,\eta
}\right\Vert _{\infty}\leq C_{\alpha}\left(  \eta\ell\left(  I\right)
\right)  ^{-\left\vert \alpha\right\vert }\left\vert I\right\vert ^{-\frac
{1}{2}}$, the rescaled function $\psi\left(  u\right)  =\left\vert
I\right\vert ^{\frac{1}{2}}h_{I}^{e,\eta}\left(  x_{I}+\ell u\right)  $ is a
smooth bump independent of $I$. Then (i) the Fresnel identity for Gaussians
gives the dispersive bound $\left\vert \int e^{-i\left(  v\cdot u+\mu
\left\vert u\right\vert ^{2}\right)  }\psi\right\vert \leq\left(  \frac{\pi
}{\left\vert \mu\right\vert }\right)  ^{\frac{n-1}{2}}\left\Vert \widehat
{\psi}\right\Vert _{L^{1}}$ uniformly in $v$, where $\widehat{\psi}\left(
\eta\right)  =\int\psi\left(  u\right)  e^{-2\pi i\eta\cdot u}du$ (write
$\psi$ by Fourier inversion and use $\left\vert \int_{\mathbb{R}}e^{-i\left(
\mu t^{2}-bt\right)  }dt\right\vert =\left(  \pi/\left\vert \mu\right\vert
\right)  ^{\frac{1}{2}}$ for $\mu\neq0$ and $b\in\mathbb{R}$), and (ii)
repeated integration by parts gives decay $C_{M}\left\vert v\right\vert ^{-M}$
when $\left\vert v\right\vert \geq4\sqrt{n-1}\left\vert \mu\right\vert $.
Together these give the radial profile%
\begin{equation}
\left\vert Th_{I}^{e,\eta}\left(  \xi\right)  \right\vert \leq C_{N}\left\vert
I\right\vert ^{\frac{1}{2}}\rho^{-\frac{n-1}{2}}\left(  1+\frac{\ell\left\vert
\xi^{\prime}+2\xi_{n}x_{I}\right\vert }{\rho}\right)  ^{-N},\ \ \ \ \ N\geq
0,\label{smooth packet}%
\end{equation}
and the proof of Proposition \ref{key} goes through with the radial kernel
$K_{r}\left(  z\right)  =r^{-\left(  n-1\right)  }\left(  1+\left\vert
z\right\vert /r\right)  ^{-2N}$, $2N>n-1$, in place of (\ref{def K}); the
center is replaced by $x\in I$ using $\left\vert x-x_{I}\right\vert \leq
\sqrt{n-1}\,r$. Lemma \ref{coef sq} holds for the Alpert coefficients
$\left\langle g,h_{I}^{e}\right\rangle $ with the same maximal function
argument. The only change is that $\triangle_{I}g$ is now a polynomial of
degree less than $\kappa$ on each child of $I$, instead of a constant. Such
functions form a finite dimensional space $V_{I}$, which is the image of
$V_{Q_{0}}$, $Q_{0}=\left[  -\frac{1}{2},\frac{1}{2}\right]  ^{n-1}$, under
the affine map taking $Q_{0}$ onto $I$, and all norms on $V_{Q_{0}}$ are
equivalent. Hence $\frac{1}{\left\vert I\right\vert }\left\Vert \varphi
\right\Vert _{L^{2}}^{2}\leq C_{n,\kappa}\left(  \frac{1}{\left\vert
I\right\vert }\int_{I}\left\vert \varphi\right\vert \right)  ^{2}$ for all
$\varphi\in V_{I}$, which replaces the computation with the children of $I$ in
the proof of Lemma \ref{coef sq}, and the Alpert square function bound of
\cite[Theorem 6]{Saw7} replaces (H2). Finally one takes $c_{I,e}=\left\langle
S_{\kappa,\eta}^{-1}f,h_{I}^{e}\right\rangle $. Conversely, the Haar case
$\kappa=1$, $\eta=0$ treated in this paper is exactly the case in which the
smoothness used in (i) and (ii) is unavailable, and the product structure of
Lemma \ref{packet} replaces it.

\subsubsection{The probabilistic estimate for the sphere\label{sphere}}

For Haar multipliers the argument above uses that the phase $x\cdot\xi
^{\prime}+\left\vert x\right\vert ^{2}\xi_{n}$ of the paraboloid is a sum of
one-variable phases. On the sphere this fails, the product profile
(\ref{packet bound}) is no longer available, and a radial profile with first
order decay is not integrable (Commentary (c) in Step 2). So the argument does
not apply directly to Haar multipliers on the sphere. However, for smooth
Alpert multipliers it does. Indeed, consider a patch of the sphere given by
$\Phi\left(  x\right)  =\left(  x,\sqrt{1-\left\vert x\right\vert ^{2}%
}\right)  $ on a small cube $S$. Write $\xi=\lambda\omega$ with $\lambda
=\left\vert \xi\right\vert $ and $\omega\in\mathbb{S}^{n-1}$. The phase
$\lambda\Phi\left(  x\right)  \cdot\omega$ is stationary in $x$ exactly when
$\omega=\pm\Phi\left(  x\right)  $, i.e.\ when $\omega$ is normal to the
sphere at $\Phi\left(  x\right)  $. The analogue of (\ref{smooth packet}) is%
\[
\left\vert Th_{I}^{e,\eta}\left(  \lambda\omega\right)  \right\vert \leq
C_{N}\left\vert I\right\vert ^{\frac{1}{2}}\rho^{-\frac{n-1}{2}}\sum_{\pm
}\left(  1+\frac{\lambda\ell\,\operatorname*{dist}\left(  \omega,\pm
\Phi\left(  x_{I}\right)  \right)  }{\rho}\right)  ^{-N},\ \ \ \ \ \rho
=1+\lambda\ell^{2}.
\]
The proof rescales $x=x_{I}+\ell u$ as before. The phase becomes $v\cdot u+\mu
Q\left(  u\right)  $, where $\mu\approx\lambda\ell^{2}$ and $Q$ has $C^{k}$
norms bounded uniformly in $\ell\leq1$, and a nondegenerate Hessian for
$\omega$ near $\pm\Phi\left(  x_{I}\right)  $. The dispersive bound is
replaced by the standard uniform decay estimate for oscillatory integrals with
nondegenerate phase \cite[Chapter VIII]{Ste2}, applied when $\left\vert
v\right\vert \lesssim\left\vert \mu\right\vert $. Precisely, what is needed is%
\[
\left\vert \int e^{-i\left(  v\cdot u+\mu Q\left(  u\right)  \right)  }%
\psi\left(  u\right)  du\right\vert \leq C\left\vert \mu\right\vert
^{-\frac{n-1}{2}},\ \ \ \ \ \left\vert \mu\right\vert \geq1,
\]
uniformly in $v\in\mathbb{R}^{n-1}$, with $C$ depending only on finitely many
derivatives of $Q$ and $\psi$ and on a lower bound for the Hessian of $Q$.
Since $Q$ is strictly convex on the support of $\psi$, the phase has at most
one critical point, and it is nondegenerate. The nonstationary bound is
unchanged in spirit: the gradient of the phase has size $\gtrsim\left\vert
v\right\vert $ when $\left\vert v\right\vert \gg\left\vert \mu\right\vert $.
Directions $\omega$ far from $\pm\Phi\left(  S\right)  $ give rapid decay. The
rest of the proof is identical, with the radial kernel of Subsection
\ref{smooth Alpert}, using polar coordinates $d\xi=\lambda^{n-1}%
d\lambda\,d\sigma\left(  \omega\right)  $ in place of $\left(  \xi^{\prime
},\xi_{n}\right)  =\left(  -2\xi_{n}z,\xi_{n}\right)  $, with $\lambda$ in
place of $\left\vert \xi_{n}\right\vert $.

\subsection{Elimination of Kakeya issues by randomization\label{kakeya}}

After Step 1, Proposition \ref{key} is the positive inequality $\left\Vert
F\right\Vert _{L^{q/2}}\lesssim\left\Vert \mathcal{S}c\right\Vert _{L^{p}}%
^{2}$, where%
\[
F=\sum_{I,e}w_{I,e}\left\vert I\right\vert \left\vert Th_{I}^{e}\right\vert
^{2},\ \ \ \ \ \left(  \mathcal{S}c\right)  ^{2}=\sum_{I,e}w_{I,e}%
\mathbf{1}_{I},\ \ \ \ \ w_{I,e}=\frac{\left\vert c_{I,e}\right\vert ^{2}%
}{\left\vert I\right\vert }\geq0.
\]
Pair $F$ with $g\geq0$ and write $\left\vert I\right\vert X=\int\mathbf{1}%
_{I}\left(  x\right)  X\,dx$:%
\[
\int_{\mathbb{R}^{n}}Fg=\int_{\mathbb{R}^{n-1}}\sum_{I,e}w_{I,e}\mathbf{1}%
_{I}\left(  x\right)  \left(  \int_{\mathbb{R}^{n}}\left\vert Th_{I}%
^{e}\right\vert ^{2}g\right)  dx\leq\left\Vert \mathcal{S}c\right\Vert
_{L^{p}}^{2}\left\Vert \mathfrak{M}g\right\Vert _{L^{\left(  p/2\right)
^{\prime}}},
\]
where%
\[
\mathfrak{M}g\left(  x\right)  \equiv\sup_{I\ni x}\max_{e\in\Gamma}%
\int_{\mathbb{R}^{n}}\left\vert Th_{I}^{e}\right\vert ^{2}g.
\]
So Proposition \ref{key} follows from the maximal inequality%
\[
\left\Vert \mathfrak{M}g\right\Vert _{L^{\left(  p/2\right)  ^{\prime}}\left(
\mathbb{R}^{n-1}\right)  }\lesssim\left\Vert g\right\Vert _{L^{\left(
q/2\right)  ^{\prime}}\left(  \mathbb{R}^{n}\right)  }.
\]
Conversely, one expects the two inequalities to be essentially equivalent by
duality, but this would require identifying the dual of the space normed by
$\left\Vert \sum_{I,e}w_{I,e}\mathbf{1}_{I}\right\Vert _{L^{p/2}}$, and we
neither need nor prove it. By Lemma \ref{packet}, $\int\left\vert Th_{I}%
^{e}\right\vert ^{2}g$ is essentially $\ell\left(  I\right)  ^{-2}$ times the
average of $g$ over the tube of $I$, plus a dispersive part. So $\mathfrak{M}$
is a maximal function over the \emph{bush} of parabolic tubes through the
origin in direction $\left(  -2x,1\right)  $. Bush estimates involve no
Kakeya-type difficulty, which is why the probabilistic problem is elementary
compared to the deterministic problem.

\end{document}